\documentclass[12pt]{article}

\usepackage{amsmath, epsfig, cite, lineno}
\usepackage{amssymb,amsthm}
\usepackage{amsfonts, color}
\usepackage{latexsym}

\newtheorem{thm}{Theorem}[section]

\newtheorem{conj}[thm]{Conjecture}
\newtheorem{lem}[thm]{Lemma}

\newcommand{\pf}{\noindent{\it Proof.} }

\numberwithin{equation}{section}

\begin{document}


\begin{center}
{\Large\bf Proof of a supercongruence conjectured by Z.-H. Sun}
\end{center}

\vskip 2mm \centerline{Victor J. W. Guo}
\begin{center}
{\footnotesize Department of Mathematics, Shanghai Key Laboratory of
PMMP, East China Normal University,\\ 500 Dongchuan Rd., Shanghai 200241,
 People's Republic of China\\
{\tt jwguo@math.ecnu.edu.cn,\quad http://math.ecnu.edu.cn/\textasciitilde{jwguo}}}
\end{center}


\vskip 0.7cm \noindent{\bf Abstract.}
The Franel numbers are defined by
$
f_n=\sum_{k=0}^n {n\choose k}^3.
$
Motivated by the recent work of Z.-W. Sun on Franel numbers, we prove that
\begin{align*}
\sum_{k=0}^{n-1}(3k+1)(-16)^{n-k-1} {2k\choose k} f_k
&\equiv 0\pmod{n{2n\choose n}}, \\
\sum_{k=0}^{p-1}\frac{3k+1}{(-16)^k} {2k\choose k} f_k
&\equiv p (-1)^{\frac{p-1}{2}} \pmod{p^3}.
\end{align*}
where $n>1$  and $p$ is an odd prime. The second congruence modulo $p^2$ confirms
a recent conjecture of Z.-H. Sun. We also show that, if $p$ is a prime of the form $4k+3$, then
$$
\sum_{k=0}^{p-1}\frac{{2k\choose k} f_k}{(-16)^k}
\equiv 0 \pmod p,
$$
which confirms a special case of another conjecture of Z.-H. Sun.



\section{Introduction}
The numbers  $f_n$, defined by
\begin{align*}
f_n=\sum_{k=0}^n {n\choose k}^3,
\end{align*}
were first studied by Franel \cite{Franel1,Franel2}, who obtained the following recurrence relation:
\begin{align*}
(n+1)^2 f_{n+1}=(7n^2+7n+2)f_{n}+8n^2 f_{n-1},\ n=1,2,\ldots. 
\end{align*}
MacMahon \cite[p.~122]{MacMahon} gave the following identity related to Franel numbers:
\begin{align}
\sum_{k=0}^{n}{n\choose k}^3x^k
=\sum_{k=0}^{n}{n+k\choose 3k}{3k\choose 2k}{2k\choose k}x^k (1+x)^{n-2k}
\label{eq:foata}
\end{align}
(see also Foata \cite{Foata} or Riordan \cite[p.~41]{Riordan}).
Strehl \cite{Strehl1} gave another expression for Franel numbers:
\begin{align*}
f_n=\sum_{k=0}^n{n\choose k}^2{2k\choose n}.
\end{align*}
Jarvis and Verrill \cite{JV} obtained the following congruence:
$$
f_{n}\equiv (-8)^n f_{p-1-n} \pmod p,
$$
where $p$ is a prime and $0\leq n\leq p-1$.
Recently, Z.-W. Sun \cite{Sun,Sun2} proved many interesting congruences involving
Franel numbers. On the other hand, during his working on congruences for Legendre polynomials \cite{SunZH1,SunZH2,SunZH3},
Z.-H. Sun \cite{SunZH00} proposed a lot of conjectures on supercongruences
concerning Franel numbers , such as
\begin{conj}\label{conj:1}{\rm(see \cite[Conjecture 4.23]{SunZH00})}Let $p>3$ be a prime. Then
\begin{align*}
\sum_{k=0}^{p-1}\frac{3k+1}{(-16)^k} {2k\choose k} f_k
\equiv p (-1)^{\frac{p-1}{2}} \pmod{p^2}.
\end{align*}
\end{conj}

\begin{conj}\label{conj:2}{\rm\cite[Conjecture 4.14]{SunZH00}}Let $p$ be an odd prime. Then
\begin{align*}
\sum_{k=0}^{p-1}\frac{{2k\choose k} f_k}{(-16)^k}
\equiv
\begin{cases}
4x^2-2p\pmod{p^2}, &\text{if $p=x^2+y^2\equiv 1\bmod{12}$ with $6\mid y$,} \\
2p-4x^2\pmod{p^2}, &\text{if $p=x^2+y^2\equiv 1\bmod{12}$ with $6\mid x-3$,} \\
4\left(\frac{xy}{3}\right)xy \pmod{p^2},
&\text{if $p=x^2+y^2\equiv 5\bmod{12}$,} \\
0 \pmod{p^2}, &\text{if $p\equiv 3\bmod 4$,}
\end{cases} 
\end{align*}
where $\left(\frac{a}{3}\right)$ is the Legendre symbol.
\end{conj}

In this paper, we shall prove the following results.
\begin{thm}\label{thm:main1} Let $n>1$ be a positive integer. Then
\begin{align}
\sum_{k=0}^{n-1}(3k+1)(-16)^{n-k-1} {2k\choose k} f_k \equiv 0\pmod{n{2n\choose n}}. \label{eq:f3}
\end{align}
\end{thm}

\begin{thm}\label{thm:main2}Let $p$ be an odd prime. Then
\begin{align}
\sum_{k=0}^{p-1}\frac{3k+1}{(-16)^k} {2k\choose k} f_k
\equiv p (-1)^{\frac{p-1}{2}} \pmod{p^3}.\label{eq:fpp}
\end{align}
\end{thm}

\begin{thm}\label{thm:main3} Let $p$ be a prime of the form $4k+3$. Then
\begin{align}
\sum_{k=0}^{p-1}\frac{{2k\choose k} f_k}{(-16)^k}
\equiv 0 \pmod p.  \label{eq:fqq}
\end{align}
\end{thm}
It is obvious that the congruence \eqref{eq:fpp} modulo $p^2$ confirms Conjecture \ref{conj:1},
while \eqref{eq:fqq} is a special case of Conjecture \ref{conj:2}.

\section{Proof of Theorem \ref{thm:main1}}
We need the following identity due to Z.-W. Sun \cite[(2.3)]{Sun}:
\begin{align}
f_n
=\sum_{k=0}^{n}{n+2k\choose 3k}{3k\choose k}{2k\choose k}(-4)^{n-k}, \label{eq:mac}
\end{align}
which can be proved by the Zeilberger algorithm (see \cite{Koepf,PWZ}). Moreover, by induction on $n$,
we can easily prove that, for all $0\leqslant k\leqslant n$,
\begin{align}
\sum_{m=k}^{n-1}(3m+1)(-16)^{n-m-1} {2m\choose m}{m+2k\choose 3k}(-4)^{m-k}
={2n\choose n}{n+2k\choose 3k}\frac{n(k-n)(-4)^{n-k}}{8(2k+1)}. \label{eq:induc}
\end{align}
In fact, when $n=k$, each side of \eqref{eq:induc} equals $0$. Suppose that \eqref{eq:induc}
is true for $n$. Then
\begin{align*}
&\hskip -2mm \sum_{m=k}^{n}(3m+1)(-16)^{n-m} {2m\choose m}{m+2k\choose 3k}(-4)^{m-k} \\
&=(3n+1){2n\choose n}{n+2k\choose 3k}(-4)^{n-k} \\
&\qquad{}+\sum_{m=k}^{n-1}(3m+1)(-16)^{n-m} {2m\choose m}{m+2k\choose 3k}(-4)^{m-k}\\
&=(3n+1){2n\choose n}{n+2k\choose 3k}(-4)^{n-k}-2{2n\choose n}{n+2k\choose 3k}\frac{n(k-n)(-4)^{n-k}}{(2k+1)} \\
&={2n+2\choose n+1}{n+2k+1\choose 3k}\frac{(n+1)(k-n-1)(-4)^{n-k+1}}{8(2k+1)}.
\end{align*}
Namely, the identity \eqref{eq:induc} holds for $n+1$.

By \eqref{eq:mac} and \eqref{eq:induc}, we have
\begin{align}
&\hskip -3mm \sum_{m=0}^{n-1}(3m+1)(-16)^{n-m-1} {2m\choose m} f_m  \notag  \\
&=\sum_{m=0}^{n-1}(3m+1)(-16)^{n-m-1} {2m\choose m} \sum_{k=0}^{m}{m+2k\choose 3k}{3k\choose k}{2k\choose k}(-4)^{m-k}.  \notag \\
&=\sum_{k=0}^{n-1}{2n\choose n}{n+2k\choose 3k}{3k\choose k}{2k\choose k}\frac{n(k-n)(-4)^{n-k}}{8(2k+1)}. \label{eq:sum1}
\end{align}
Note that, for $n\geqslant 2$ and $0\leqslant k<n$, both
$$
\frac{1}{2k+1}{3k\choose k}={3k\choose k}-2{3k\choose k-1}$$
and
$
{2k\choose k}\frac{(-4)^{n-k}}{8}
$
are integers. From \eqref{eq:sum1} we deduce that
\begin{align}
\frac{1}{n{2n\choose n}}\sum_{k=0}^{n-1}(3k+1)(-16)^{n-k-1} {2k\choose k} f_k
&=\sum_{k=0}^{n-1}{n+2k\choose 3k}{3k\choose k}{2k\choose k}\frac{(k-n)(-4)^{n-k}}{8(2k+1)}\label{eq:newsum}
\end{align}
is an integer. This proves \eqref{eq:f3}.

\section{Proof of Theorem \ref{thm:main2}}
It follows from
\eqref{eq:newsum} that
\begin{align}
\frac{1}{p}\sum_{k=0}^{p-1}\frac{3k+1}{(-16)^k} {2k\choose k} f_k
&=-\,4^{1-p}{2p-1\choose p-1}\sum_{k=0}^{p-1}{p+2k\choose 3k}{3k\choose k}{2k\choose k}\frac{k-p}{(2k+1)(-4)^k}. \label{eq:newsum-pp}
\end{align}
Note that Babbage \cite{Babbage} proved the following congruence:
\begin{align}
{2p-1\choose p-1}&\equiv 1 \pmod{p^2}. \label{eq:Babbage}
\end{align}
Moreover, we have
\begin{align}
{p+2k\choose 3k}{3k\choose k}&=\frac{(p+2k)!}{(2k)! k!(p-k)!}
\equiv
\begin{cases}(-1)^{k-1} \frac{p}{k}, &\text{if $1\leqslant k<\frac{p-1}{2},$} \\[5pt]
(-1)^{k-1} \frac{2p}{k}, &\text{if $\frac{p+1}{2}\leqslant k<p$}
\end{cases} \pmod{p^2},  \label{eq:222}
\end{align}
and, when $k=\frac{p-1}{2}$,
\begin{align}
{p+2k\choose 3k}{3k\choose k}{2k\choose k}\frac{k-p}{2k+1}=-{2p-1\choose p-1}{p-1\choose \frac{p-1}{2}}^2
\equiv -16^{p-1} \pmod{p^2},    \label{eq:333}
\end{align}
where we have used Babbage's congruence \eqref{eq:Babbage} and Morley's congruence \cite{Morley}:
\begin{align*}
{p-1\choose \frac{p-1}{2}}\equiv (-1)^{\frac{p-1}{2}} 4^{p-1} \pmod{p^3}.
\end{align*}
Substituting \eqref{eq:Babbage}--\eqref{eq:333} into \eqref{eq:newsum-pp},
and observing that ${2k\choose k}\equiv 0\pmod{p}$ for $\frac{p-1}{2}<k<p$, we obtain
\begin{align}
\frac{1}{p}\sum_{k=0}^{p-1}\frac{3k+1}{(-16)^k} {2k\choose k} f_k
&\equiv p4^{1-p}+(-4)^{\frac{p-1}{2}}+4^{1-p}\sum_{k=1}^{\frac{p-3}{2}}{2k\choose k}\frac{p-p^2/k}{(2k+1)4^k}  \nonumber \\
&\equiv (-4)^{\frac{p-1}{2}}+4^{1-p}\sum_{k=0}^{\frac{p-3}{2}}{2k\choose k}\frac{p}{(2k+1)4^k}  \pmod{p^2}. \label{eq:newsum2}
\end{align}
Since
\begin{align}
{2k\choose k}4^{-k}\equiv (-1)^k {\frac{p-1}{2} \choose k} \pmod{p}, \label{eq:pmod}
\end{align}
we may rewrite \eqref{eq:newsum2} as
\begin{align}
\frac{1}{p}\sum_{k=0}^{p-1}\frac{3k+1}{(-16)^k} {2k\choose k} f_k
&\equiv (-4)^{\frac{p-1}{2}}+4^{1-p}\sum_{k=0}^{\frac{p-3}{2}} (-1)^k {\frac{p-1}{2}\choose k}\frac{p}{2k+1} \pmod{p^2},  \label{eq:newsum3}
\end{align}
Applying the famous identity
$$
\sum_{k=0}^{n}(-1)^k{n\choose k}\frac{1}{x+k}=\frac{n!}{x(x+1)\cdots (x+n)}
$$
with $x=\frac{1}{2}$ and $n=\frac{p-1}{2}$, we may simplify \eqref{eq:newsum3} as
\begin{align*}
\frac{1}{p}\sum_{k=0}^{p-1}\frac{3k+1}{(-16)^k} {2k\choose k} f_k
&\equiv (-4)^{\frac{p-1}{2}}+2^{1-p}{p-1\choose \frac{p-1}{2}}^{-1}-(-1)^{\frac{p-1}{2}}4^{1-p} \\
&\equiv (-1)^{\frac{p-1}{2}}\left(2^{p-1}+8^{1-p}-4^{1-p}\right) \pmod{p^2}.
\end{align*}
By Fermat's little theorem, we have $2^{p-1}-1\equiv 0\pmod p$ and so
$$
2^{p-1}+8^{1-p}-4^{1-p}-1=(2^{p-1}-1)^2(2^{1-p}+4^{1-p}+8^{1-p})\equiv 0 \pmod{p^2}.
$$
This completes the proof. \qed

\section{Proof of Theorem \ref{thm:main3}}
We first give a binomial coefficient identity.
\begin{lem}Let $n$ and $k$ be nonnegative integers with $k\leqslant n$. Then
\begin{align}
\sum_{m=k}^{n}{n\choose m}{m+2k\choose 3k}(-1)^{m-k}={2k\choose n-k}(-1)^{n-k}. \label{eq:sum-new}
\end{align}
\end{lem}
\pf Consider the generating function of the left-hand side of \eqref{eq:sum-new}. By the binomial theorem, we have
\begin{align*}
\sum_{n=k}^\infty x^n \sum_{m=k}^{n}{n\choose m}{m+2k\choose 3k}(-1)^{m-k}
&=\sum_{m=k}^\infty {m+2k\choose 3k}(-1)^{m-k} \sum_{n=m}^{\infty}{n\choose m} x^n \\
&=\sum_{m=k}^\infty {m+2k\choose 3k} \frac{(-1)^{m-k} x^m}{(1-x)^{m+1}}  \\
&=\frac{x^k}{(1-x)^{x+1}}\left(1+\frac{x}{1-x}\right)^{-3k-1} \\
&=x^k(1-x)^{2k},
\end{align*}
which is clearly the generating function of the right-hand side of \eqref{eq:sum-new}.
\qed

By \eqref{eq:pmod} and \eqref{eq:mac}, we have
\begin{align}
\sum_{k=0}^{p-1}\frac{{2k\choose k} f_k}{(-16)^k}
&\equiv \sum_{m=0}^{p-1}\frac{{\frac{p-1}{2}\choose m} f_m}{4^m}  \nonumber \\
&=\sum_{m=0}^{p-1} \frac{{\frac{p-1}{2}\choose m}}{4^m} \sum_{k=0}^{m}{m+2k\choose 3k}{3k\choose k}{2k\choose k}(-4)^{m-k}
\pmod p.  \label{eq:doub-sum}
\end{align}
Exchanging the summation order in the right-hand side of \eqref{eq:doub-sum}, and applying \eqref{eq:sum-new} with $n=\frac{p-1}{2}$,
we obtain
\begin{align}
\sum_{k=0}^{p-1}\frac{{2k\choose k} f_k}{(-16)^k}
\equiv
\sum_{k=0}^{\frac{p-1}{2}}\frac{(-1)^{\frac{p-1}{2}-k}}{4^k}{3k\choose k}{2k\choose k}{2k\choose {\frac{p-1}{2}}-k} \label{eq:final-0}
\pmod p.
\end{align}
It is easy to see that
\begin{align}
{2k\choose {\frac{p-1}{2}}-k}\equiv {\frac{3(p-1)}{2}-3k\choose {\frac{p-1}{2}}-k} (-1)^{\frac{p-1}{2}-k} \pmod p. \label{eq:final-2}
\end{align}
Thus, by \eqref{eq:pmod} and \eqref{eq:final-2}, we may rewrite \eqref{eq:final-0} as
\begin{align}
\sum_{k=0}^{p-1}\frac{{2k\choose k} f_k}{(-16)^k}
\equiv
\sum_{k=0}^{\frac{p-1}{2}} (-1)^k {\frac{p-1}{2}\choose k}{3k\choose k}{\frac{3(p-1)}{2}-3k\choose {\frac{p-1}{2}}-k}
\pmod p.  \label{eq:final-3}
\end{align}
If $p$ is a prime of the form $4k+3$, then $\frac{p-1}{2}$ is odd, and by the symmetry of binomial coefficients, the right-hand side of \eqref{eq:final-3} equals $0$.
This completes the proof.

\section{Concluding remarks and open problems}
Inspired by Conjectures 4.22 and 4.23 in \cite{SunZH00} and Theorem \ref{thm:main1} in this paper, we
propose the following two conjectures.
\begin{conj}\label{conj:new1}Let $n>1$ be a positive integer and let
$(a,b,c)\in\{(9,4,5),(5,2,16),\break(9,2,50), (5,1,96),(6,1,320),(90,13,896),(102,11,10400)\}$.
Then
\begin{align*}
\sum_{k=0}^{n-1}(ak+b)c^{n-k-1} {2k\choose k} f_k \equiv 0\pmod{n{2n\choose n}}.
\end{align*}
\end{conj}

\begin{conj}\label{conj:new2}Let $n>1$ be a positive integer and let
$(a,b,c)\in\{(15,4,-49),(9,2,-112),\break(99,17,-400),(855,109,-2704),(585,58,-24304)\}$.
Then
\begin{align*}
\sum_{k=0}^{n-1}(ak+b)c^{n-k-1} {2k\choose k} f_k \equiv 0\pmod{n{2n\choose n}}.
\end{align*}
\end{conj}
We have verified Conjectures \ref{conj:new1} and \ref{conj:new2} for $n$ up to $500$ via Maple.

We also make the third conjecture on Franel numbers as follows.
\begin{conj} Let $m,n$ be positive integers and let $a_1,\ldots,a_m$ be integers.
Then
\begin{align}
\sum_{k=0}^{n-1}(3k+2)(-1)^{(m-1)k} f_k \prod_{i=1}^{m}{a_i n-1\choose k}{a_i n+k\choose k}&
\equiv 0 \pmod{n^2}, \label{eq:conj-1}\\
\sum_{k=0}^{n-1}(9k^2+5k)(-1)^{(m-1)k} f_k \prod_{i=1}^{m}{a_i n-1\choose k}{a_i n+k\choose k}&
\equiv 0 \pmod{n^2}.  \label{eq:conj-2}
\end{align}
\end{conj}
Note that, for any prime $p$ and nonnegative integer $k\leqslant p-1$, there holds
\begin{align*}
{a_i p-1\choose k}{a_i p+k\choose k}=\prod_{j=1}^{k}\frac{(a_i p-j)(a_i p+j)}{j^2}\equiv (-1)^k \pmod{p^2}.
\end{align*}
Therefore, by \cite[Theorem 1.1]{Sun}, the congruences \eqref{eq:conj-1} and \eqref{eq:conj-2} are true for any prime $n$.
On the other hand, the author \cite{Guo} has  proved that
\begin{align*}
\sum_{k=0}^{n-1}(3k+2)(-1)^k f_k \equiv 0 \pmod{2n^2},
\end{align*}
which was conjectured by Z.-W. Sun (see \cite[Conjecture 1.3]{Sun2}); and for $m=0$ and $n>1$, the congruence \eqref{eq:conj-2}
can be further strengthened as
\begin{align*}
\sum_{k=0}^{n-1}(9k^2+5k)(-1)^k f_k \equiv 0 \pmod{n^2(n-1)},
\end{align*}
which was also conjectured by Z.-W. Sun (see Conjecture 5.3 in a previous version of \cite{Sun2}: {\tt http://arxiv.org/pdf/1112.1034v10.pdf}).


\vskip 5mm
\noindent{\bf Acknowledgments.} This work was partially
supported by the Fundamental Research Funds for the Central Universities and
the National Natural Science Foundation of China (grant 11371144).


\begin{thebibliography}{99}
\small \setlength{\itemsep}{-.8mm}

\bibitem{Babbage}C. Babbage, Demonstration of a theorem relating to prime numbers,
Edinburgh Philos. J. 1 (1819), 46--49.

\bibitem{Foata}D. Foata, Etude alg\'ebrique de certains probl\`emes d'analyse combinatoire et du
calcul des probabilit\'es, Publ. Inst. Statist. Univ. Paris 14 (1965), 81--241.

\bibitem{Franel1}J. Franel, On a question of Laisant, L'interm\'ediaire des math\'ematiciens 1 (1894), 45--47.

\bibitem{Franel2}J. Franel, On a question of J. Franel, L'interm\'ediaire des math\'ematiciens 2 (1895), 33--35.

\bibitem{Guo}V.J.W. Guo, Proof of two conjectures of Sun on congruences for Franel numbers,
Integral Transforms Spec. Funct. 24 (2013), 532--539.

\bibitem{JV}F. Jarvis  and H. A. Verrill, Supercongruences for the Catalan-Larcombe-French numbers,
Ramanujan J. 22 (2010), 171--186.

\bibitem{Koepf}W. Koepf, Hypergeometric Summation, an Algorithmic Approach to Summation
and Special Function Identities, Friedr. Vieweg \& Sohn, Braunschweig, 1998.

\bibitem{MacMahon}P. A. MacMachon, Combinatorial Analysis, Vol. 1, Cambridge University Press, London, 1915.

\bibitem{Morley}F. Morley, Note on the congruence $2^{4n}\equiv (-1)^n (2n)!/(n!)^2$, where $2n+1$ is a prime,
Ann. Math.  9 (1894/95), 168--170.

\bibitem{PWZ}M. Petkov\v{s}ek, H. S. Wilf, and D. Zeilberger, $A=B$, A K Peters, Ltd., Wellesley, MA, 1996.

\bibitem{Riordan}J. Riordan, Combinatorial Identities, J. Wiley, New York, 1979.

\bibitem{Strehl1}V. Strehl, Binomial sums and identities, Maple Technical Newsletter 10 (1993), 37--49.


\bibitem{SunZH1}Z.-H. Sun, Congruences concerning Legendre polynomials, Proc. Amer. Math. Soc. 139 (2011), 1915--1929.

\bibitem{SunZH2}Z.-H. Sun, Congruences concerning Legendre polynomials II, J. Number Theory 133 (2013), 1950--1976.

\bibitem{SunZH3}Z.-H. Sun, Congruences concerning Legendre polynomials III, Int. J. Number Theory  9 (2013), 965--999.

\bibitem{SunZH00}Z.-H. Sun, Some conjectures on congruences, preprint, 2013, arXiv:1103.5384v5.


\bibitem{Sun}Z.-W. Sun, Connections between $p=x^2+3y^2$ and Franel numbers, J. Number Theory 133 (2013), 2914--2928.

\bibitem{Sun2}Z.-W. Sun, Congruences for Franel numbers, Adv. Appl. Math. 51 (2013), 524--535.


\end{thebibliography}
\end{document}